\documentclass[twoside,12pt, leqno]{article}
\usepackage{amsmath,amscd,amsthm,amssymb,amsxtra,latexsym,epsfig,epic,graphics}
\usepackage[matrix,arrow,curve]{xy}
\usepackage{graphicx}
\usepackage{diagrams}
\def\antiddot{\mathinner{\mkern1mu\raise1pt\vbox{\kern7pt\hbox{.}}\mkern2mu
        \raise4pt\hbox{.}\mkern2mu\raise7pt\hbox{.}\mkern1mu}}

\newcommand{\KK}{{\mathbb K}}
\newcommand{\PP}{{\mathbb P}}

\newcommand{\TT}{{\mathbb T}}

\newcommand{\ZZ}{{\mathbb Z}}

\newcommand{\ann}{{\rm{ann}}}
\newcommand{\Ext}{{\rm{Ext}}}

\newcommand{\s}{\mathcal}

\newcommand{\sB}{{\s B}}

\newcommand{\sF}{{\s F}}
\newcommand{\sG}{{\s G}}

\newcommand{\sM}{{\s M}}

\newcommand{\sO}{{\s O}}

\newcommand{\sV}{{\s V}}

\newcommand{\sZ}{{\s Z}}

\newcommand{\tensor}{\otimes}

\newcommand{\punkt}{\hspace{-.3ex}\raise.15ex\hbox to1ex{\Huge.}}

\DeclareMathOperator{\reg}{reg}

\DeclareMathOperator{\Hom}{Hom}

\newtheorem{theorem}{Theorem}[section]
\newtheorem{lemma}[theorem]{Lemma}
\newtheorem{proposition}[theorem]{Proposition}
\newtheorem{corollary}[theorem]{Corollary}

\theoremstyle{definition}

\newtheorem{example}[theorem]{Example}

\def\cF{{\cal F}}

\def\BS{{Boij-S\"oderberg\ }}
\DeclareMathOperator{\rH}{{\rm H}}
\def\bC{{\mathbb C}}
\def\Gr{{\rm Gr}}

\makeatletter
\def\Ddots{\mathinner{\mkern1mu\raise\p@
\vbox{\kern7\p@\hbox{.}}\mkern2mu
\raise4\p@\hbox{.}\mkern2mu\raise7\p@\hbox{.}\mkern1mu}}
\makeatother

\DeclareMathOperator{\coreg}{coreg}

\makeatletter
\def\Ddots{\mathinner{\mkern1mu\raise\p@
\vbox{\kern7\p@\hbox{.}}\mkern2mu
\raise4\p@\hbox{.}\mkern2mu\raise7\p@\hbox{.}\mkern1mu}}
\makeatother

\usepackage{times}
\newdimen\x \x=12pt

\usepackage{color}

\date{}
\title{The Banks of the Cohomology River}
\author{David Eisenbud and Frank-Olaf Schreyer
\footnote{This paper reports on work done during a period
of ``Research in Pairs'' at the Mathematische Forschungsinstitute
Oberwolfach, August 1--14, 2011. We are grateful
to Institute for providing such a beautiful and peaceful setting, together
with wonderful resources for research. The first author is also grateful to the 
National Science Foundation for partial support during this period.}}

\begin{document}

\maketitle

{\it To the Memory of Masaki Maruyama}
\begin{abstract}
We  give sharp bounds on the vanishing of the cohomology of a tensor product of vector bundles on $\PP^{n}$ in terms of the vanishing of the cohomology of the factors. For this purpose we introduce regularity indices generalizing the Castelnuovo-Mumford regularity.

As an application we give a sufficient condition for a vector bundle to have an unobstructed deformation theory that depends only on the cohomology table
of the bundle. We construct complete families
of bundles with such cohomology tables.
\end{abstract}

\section{Introduction}

If we know which cohomology groups of (all twists of) 
two vector bundles $\sF, \sG$ on $\PP^{n}$ are zero and nonzero,
what can we say about the cohomology of twists of $F\otimes G$?  For example, one might naively suppose that if $H^{i}\sF(a) \neq 0$
and $H^{j}\sG(b) \neq 0$, and if $i+j\leq n$, then at least for some sheaves with the given vanishing pattern one might have $H^{i+j}\sF\otimes \sG(a+b)\neq 0$. 

In this paper we will give sharp bounds on which cohomology
groups of twists of $\sF\otimes \sG$ vanish, and we will see that they are much more
restrictive than  the naive idea above would suggest (see Example~\ref{tensor product example}). 

We were led to these bounds by a result from the \BS theory of cohomology tables of vector bundles: by Theorems 0.5 and 6.2 of our paper \cite{ES1}
there is, for any vector bundle $\sF$ on projective space a uniquely defined homogeneous vector bundle (that is, a direct sum of twists of Schur functors applied to
the tangent bundle) that has, in characteristic zero, the same cohomology table as
$\sF$ up to a rational multiple.
The inequalities we prove are sharp for these homogeneous bundles in characteristic zero, in a rather strong sense.

The inequalities we give strengthen those of Sidman~\cite{sidman} and 
Caviglia~\cite{caviglia}. Those authors' work is based on the characterization of regularity in terms of ``approximate'' free resolutions---that is, free complexes that are  resolutions away from some low-dimensional locus (the idea of using such approximate resolutions seems to go back to the paper of Gruson-Lazarsfeld-Peskine~\cite{GLP}.) The  improvement that gives us stronger bounds is the use of (approximate) free monads instead of  resolutions. These ideas are described in \S~\ref{sec:monads}. 

One of the most interesting tensor products of two bundles is 
$End(\sF) = \sF^{*} \otimes \sF$, and one of its most interesting cohomology groups is 
$$
H^{2}(\sF^{*} \otimes \sF) = \Ext^{2}(\sF, \sF),
$$
the obstruction space for deformations of $\sF$. The bounds on cohomology of a tensor product allow us to give an interesting sufficient condition under which this deformation space is zero, so that the local deformation space of $\sF$ is smooth. It turns out that in these unobstructed cases we can actually write down complete families of the bundles with the given cohomology that are irreducible smooth rational varieties. This application occupies \S~\ref{sec:natural families}.

The bounds will be given in terms of \emph{regularity indices}, which we will now describe.

\subsection{Regularity Indices}

Let $\KK$ be a field, and let $\sF$ be a coherent sheaf on $\PP^n=\PP^n_\KK$.
We define
 the \emph{$k$-th regularity index} of $\sF$ to be
$$
\reg^k\sF := \inf \{ m \mid H^{j} \sF(m-j) = 0 \hbox{ for all } j >k \}.
$$
and the $k$-th coregularity index to be
$$
\coreg^{k}\sF := \sup \{ m \mid H^{j} \sF(m-j) = 0 \hbox{ for all } j < n-k  \}.
$$
(This definition \emph{differs} from the one in \cite{ES1}!) Thus $\reg^0 \sF \le m$ if and only if $\sF$ is $m$-regular in the classical sense of Castelnuovo-Mumford.
For any vector bundle $\sF$ and any integer $m$ we have $\coreg^{m}\sF = -\reg^{m}\sF^{*} - 1$, as one
sees easily by duality.

The \emph{cohomology table} of $\cF$ is the collection of numbers
$$ 
\gamma(\cF) =\{h^{i}(\sF(d)): = \dim_{\KK} H^i(\cF(d)\}.
$$
We display it in a table with $h^{i}(\sF(d))$ in the $i$-th row (numbering from the bottom) and the $(i+d)$-th column (numbering from left to right), and to simplify the picture we replace the zeros elements by dots. As explained in our \cite{EFS} the
cohomology table, in this form, is also the Betti table of the Tate resolution associated
to $\sF$, and it follows that
If $\reg^{k}\sF = m$ then (as in the case of the Castelnuovo-Mumford regularity,
$k=0$) we have 
$H^{j} \sF(m'-j) = 0 \hbox{ for all } j >k$ for every $m'> m$, and similarly for the 
coregularity (Proposition~\ref{inductive regularity}.)

For example, the cohomology table  of the Horrocks-Mumford bundle $\sF_{HM}$ on $\PP^{4}$ is:
\goodbreak
\begin{verbatim}
4: 100 35  4  .  . . .  .  .  .   .
3:   .  2 10 10  5 . .  .  .  .   .
2:   .  .  .  .  . 2 .  .  .  .   .
1:   .  .  .  .  . . 5 10 10  2   .
0:   .  .  .  .  . . .  .  4 35 100
    -5 -4 -3 -2 -1 0 1  2  3  4   5
\end{verbatim}
In this display, $\reg^{m}\sF$ is the number of the
 left-most column with only dots above row
$m$, and $\coreg^{m}\sF$ is the number of the right-most column
with only dots below row $n-m$. Thus for example $\reg^{1}\sF_{HM}=1$ and $\coreg^{0}\sF_{HM} = -5$

\subsection{Banks of the Cohomology River}
We think of the nonzero values of the cohomology of $\sF$ as forming the
\emph{cohomology river}, and the regularities and coregularities as defining its \emph{banks}. Our first main result describes the banks of the cohomology river of a tensor product:

\begin{theorem}\label{main2-intro}If $\sF$ and $\sG$ are vector bundles on $\PP^n$ then
$$
\reg^p(\sF\tensor \sG) \le \min_{ k+l=p} \bigl(\reg^k\sF +\reg^\ell \sG \bigr)
$$
and 
$$
\coreg^p(\sF\tensor \sG) \geq 1+ \max_{ k+l=p} \bigl(\coreg^k\sF +\coreg^\ell \sG \bigr). 
$$
The  inequality  for $\reg^{p}(\sF\otimes \sG)$ holds, more generally, for any
coherent sheaves $\sF$ and $\sG$ such that
the dimension of the sheaf $Tor_{1}(\sF,\sG)$ is at most $p+2$.
\end{theorem}

Our second main result shows that Theorem~\ref{main2-intro} is sharp in a strong sense:
\begin{theorem}\label{sharpness}
Given any pair of cohomology tables $\Phi, \Gamma$ of vector bundles on $\PP^{n}$, there exists a pair of homogeneous vector bundles
$\sF$ and $\sG$ on $\PP^{n}_{\mathbb C}$ whose cohomology tables are rational multiples of $\Phi$ and $\Gamma$, and such that equality holds for every $p$ in the formulas for
$\reg^{p}(\sF\otimes \sG)$ and $\coreg^{p}(\sF\otimes \sG)$ 
of Theorem~\ref{main2-intro}. 
\end{theorem}
The proofs are given in \S~\ref{sec:main proofs}.
\begin{example}\label{tensor product example} In fact, the formulas of Theorem~\ref{main2-intro} seem to be
sharp rather often. For example, let $\pi: \PP^{1}\times \PP^{1}\times \PP^{1} \to \PP^{3}$ be
the projection defined by the symmetric functions as in \cite{ES1}. 
Taking $\sF = \pi_{*}\sO(4,1,-1)$ and $\sG = \pi_{*}\sO(3,-1,-2)$ we get 
bundles with cohomology tables:
\begin{verbatim}
3: 70 24  .  . .  .  .   .
2:  .  .  8  6 .  .  .   .
1:  .  .  .  . 4  .  .   .
0:  .  .  .  . . 18 56 120
   -4 -3 -2 -1 0  1  2   3
\end{verbatim}
and
\begin{verbatim}
3: 168 84 30  .  . .  .  .
2:   .  .  . 12 12 6  .  .
1:   .  .  .  .  . .  .  .
0:   .  .  .  .  . . 12 42
    -4 -3 -2 -1  0 1  2  3.
\end{verbatim}
Computing the cohomology table of the tensor product in Macaulay2 (\cite{M2}) we get:
\begin{verbatim}
3: 624 216  72   8  .  .  .   .
2:   .  96 140 144 96 42  .   .
1:   .   .   .   . 18 36 48   .
0:   .   .   .   .  .  . 24 216
    -4  -3  -2  -1  0  1  2   3
\end{verbatim}
Inspection shows that equality is achieved here for all the bounds of Theorem~\ref{main2-intro}. Notice that we have (for example) $H^{1}\sF(-1)=4$ and
$H^{2}\sG(-1) = 6$, but that no sheaves with these vanishing patterns can have
$H^{3}(\sF\otimes \sG(-2))\neq 0.$
\end{example}

%\subsection{Deformations and Obstructions}
Two values of the cohomology of a tensor product of bundles are particularly interesting: $H^{1}(\sF^{*}\otimes \sF) = \Ext^{1}(\sF, \sF)$ is the space of
first order deformations of $\sF$, and $H^{2}(\sF^{*}\otimes \sF) = \Ext^{2}(\sF, \sF)$
is the obstruction space. 
 
As an application of our theory we determine, in \S~\ref{sec:natural families}, the cohomology tables
of bundles that force the obstruction space to be trivial, and we compute rational families
of  bundles containing all bundles with such cohomology tables.

\section{\BS Theory for Vector Bundles}\label{sec:BS} 
By a  \emph{homogeneous bundle} on $\PP^{n}$ we mean the result
of applying a Schur functor $S_{\lambda}$ to the universal $n$-quotient
bundle $Q$, and then (possibly) tensoring with a line bundle. 
Here
$\lambda = \lambda_{n-1},\dots, \lambda_{0}$
is a partition with $n$ parts; that is, the
$\lambda_{i}$ are integers such that $ \lambda_{n-1}\geq \dots \geq \lambda_{0}\geq 0$. 
We choose
our conventions so that
$S_{m,0,\dots,0}Q$ is the $m$-th symmetric power of $Q$ while
$S_{1^{m}, 0, \dots,0}Q$ is the $m$-th exterior power of $Q$.
We draw the  Young diagram corresponding to
$\lambda$ by putting $\lambda_{i}$ boxes in the $i$-th row and right-justifiying the picture; for example the partition $(7,5,2,2,0,0)$
corresponds to the diagram
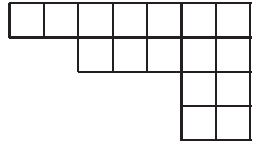
\smallbreak
\centerline{
\vbox{
\def\star{\rlap{\hbox to 13pt{\hfil\raise3.5pt\hbox{$*$}\hfil}}}
\def\ {\hbox to 13pt{\vbox to 13pt{}\hfil}}
\def\*{\star\ }
\def\_{\hbox to 13pt{\hskip-.2pt\vrule\hss\vbox to 13pt{\vskip-.2pt
            \hrule width 13.4pt\vfill\hrule\vskip-.2pt}\hss\vrule\hskip-.2pt}}
\def\x{\star\_}
\offinterlineskip
\hbox{\ \ \   \ \_\_\_\_\_\_\_\ \ \ }
\hbox{\ \ \   \  \ \ \_\_\_\_\_\ \ \ }
\hbox{\ \ \   \  \ \ \ \ \  \_\_\ \ \ }
\hbox{\ \ \   \  \  \  \  \  \  \_\_\ \ \ }
\hbox{\ \ \   \  \  \  \  \  \  \  \  \ \ \ }
}
}  
\noindent 
(where rows 0 and 1 have zero boxes!)

The vanishing part of Bott's Theorem about homogeneous bundles in characteristic zero has a very simple statement in terms of cohomology tables (this was pointed out to us by Jerzy Weyman). 

\begin{theorem} [Bott]\label{Bott} Let $\lambda_{n-1},\dots,\lambda_{0}$ be a partition as above, and let $Q$ be the universal rank $n$ quotient
bundle on $\PP^{n}_{\mathbb C}$. The cohomology table of $S_{\lambda}(Q)$ has the form

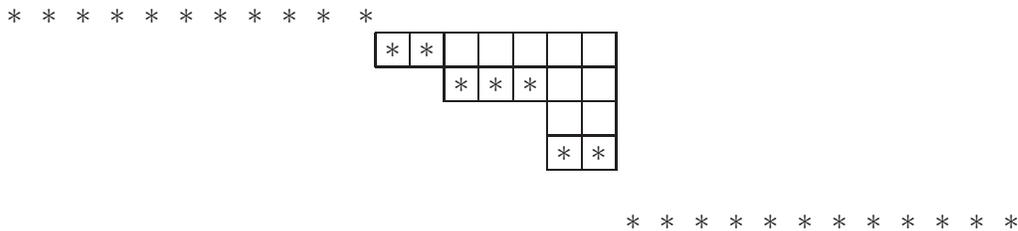
\centerline{
\vbox{
\def\star{\rlap{\hbox to 13pt{\hfil\raise3.5pt\hbox{$*$}\hfil}}}
\def\ {\hbox to 13pt{\vbox to 13pt{}\hfil}}
\def\*{\star\ }
\def\_{\hbox to 13pt{\hskip-.2pt\vrule\hss\vbox to 13pt{\vskip-.2pt
            \hrule width 13.4pt\vfill\hrule\vskip-.2pt}\hss\vrule\hskip-.2pt}}
\def\x{\star\_}
\offinterlineskip
\hbox{\*\*\*\*\*\*\*\*\*\* \* \ \ \ \ \ \ \ \ \ \ \ }
\hbox{\ \ \ \ \ \  \\ \ \   \\ \ \   \ \x\x\_\_\_\_\_\ \ \ \ \ \ \ \ \ \ \ \ }
\hbox{\ \ \ \   \\ \ \   \\ \ \   \\ \ \   \  \ \ \x\x\x\_\_\ \ \ \ \ \ \ \ \ \ \ \ }
\hbox{\ \ \ \   \\ \ \   \\ \ \   \\ \ \   \  \ \ \ \ \  \_\_\ \ \ \ \ \ \ \ \ \ \ \  }
\hbox{\ \ \ \   \\ \ \   \\ \ \   \\ \ \   \  \  \  \  \  \  \x\x\ \ \ \ \ \ \ \ \ \ \ \ }
\hbox{\ \ \ \   \\ \ \   \\ \ \   \\ \ \   \  \  \  \  \  \  \  \  \ \ \ \ \ \ \ \ \ \ \ \ }
\hbox{\ \ \ \   \\ \ \   \\ \ \   \\ \ \   \   \  \  \  \  \  \  \  \*\*\*\*\*\*\*\*\*\*\*\* }
}
}  
\noindent 
where the nonzero entries of the table are exactly those marked by $*$,
the top row of the Young diagram is row $n-1$,
and the right-hand column of the Young diagram is  column $-1$.
\end{theorem}

For example, we see from Theorem~\ref{Bott}  that 
$$
\reg^{k} S_{\lambda}(Q) = -\lambda_{k}.
$$

We partially order the partitions component-wise---in terms of Young diagrams this is the partial order by inclusion. One of the main results of \BS theory for vector bundles
can be thought of as associating to any vector bundle on projective space
a homogeneous bundle with (in characteristic zero) the same cohomology table, up
to a rational multiple. We restrict ourselves to 0-regular bundles for simplicity; of course we can
apply the result to any bundle by first tensoring with a sufficiently positive line bundle. The following statement combines Theorems 0.5 and 6.2 of our paper \cite{ES1}.

\begin{theorem}\label{BS for VB}
The cohomology table of any bundle $\sF$ with $\reg^0 \sF \le 0$ can  be written uniquely as a positive rational linear combination of the (characteristic zero) cohomology tables of a sequence of homogeneous bundles corresponding to a totally ordered set of Young diagrams. 
\end{theorem}

\noindent{\bf Remark:} One can use the \BS decomposition to bound the numbers in the cohomology table of the tensor product using just the knowledge
of which entries of the cohomology table are zero and the Hilbert polynomial.
But one might hope for a still stronger principle, asserting perhaps that
if the cohomology tables $\Phi$ and $\Gamma$
of bundles $F$ and $G$ have \BS decompositions
$$
\Phi = \sum \alpha_{i}\Phi^{i},\quad \Gamma =\sum\beta_{j}\Gamma^{j},
$$
where $\Phi^{i}$ and $\Gamma^{j}$ are the cohomology tables of
the homogeneous bundles $S_{\phi^{i}}Q$ and $S_{\gamma^{j}}Q$, then
the cohomology table of $F\otimes G$ would be bounded, term by term, by the sum
of $\alpha_{i}\beta_{j}$ times the  cohomology table  of
$S_{\phi^{i}}Q\otimes S_{\gamma^{j}}Q$. This is false, as
Example~\ref{tensor product example} shows. 

\section{Linear Monads and  Regularity Indices}\label{sec:monads}

It is well known that the (Castelnuovo-Mumford) regularity 
$\reg^{0}\sF$ of a coherent sheaf $\sF$ on $\PP^{n}$ can be characterized as the smallest integer $m$
such that $\sF(m)$ admits a \emph{linear free resolution}, that is, such that there
is a complex
$$
\sM:  \cdots\to \sO_{\PP^n}(-t-1)^{\beta_{t+1}}\to \sO_{\PP^n}(-t)^{\beta_{t}} \to\cdots\to \sO_{\PP^n}^{\beta_{0}}\to 0
$$
with homology $ \rH^{0}(\sM) = \sF$ and no other homology.
Our next result is a characterization of this sort for all the regularity indices of a sheaf and also the coregularity indices of a vector bundle. This is the characterization that we will use to prove our main theorem. 

Recall that a \emph{monad} $\sM$ for a sheaf $\sF$ is a finite complex of sheaves
$$ 
\cdots \to \sM^{-1} \to \sM^{0} \to \sM^{1} \to \dots 
$$
whose only homology is
$\rH^*(\sM) = \rH^0(\sM) \cong \sF$. The monad is called \emph{linear} if
$\sM^{i}$ is a direct sum of copies of $\sO(i)$ for each $i$.

\begin{proposition}\label{reg from monad} If $\sF$ is a coherent sheaf on $\PP^{n}$ then
$\reg^{k}\sF$ is the smallest integer $m$ such that $\sF(m)$ admits a linear monad
$\sM$ with $\sM^{\ell}=0$ for all $\ell>k$. 
If $\sF$ is a vector bundle then $\coreg^{k}\sF$ is the largest integer $m$ such that $\sF(m+ 1)$ admits a linear
monad $\sM$ with $\sM^{\ell}=0$ for all $\ell<-k$. 
\end{proposition}

\begin{proof} 
If $\sF$ is a vector bundle then the dual of a monad for $\sF$ is a monad for $\sF^{*}$. Using the formula $\coreg^{k}\sF = -\reg^{k} \sF^{*}-1$ we see that the second statement of the Proposition follows from the first.
%%\frank{Is it true that the second part holds without the assumption that $\sF$ is a vector bundle?
%I do not see how Lemma 2.2 could apply. The chase ends up with an uncontroled $H^0$
%}

Twisting by $-m$, the first statement will follow if we show that
a coherent sheaf $\sF$ admits a linear monad $\sM$ with $\sM^{\ell} = 0$
for all $\ell>k$ if and only if $\reg^{k}\sF\leq 0$. The ``only if'' part follows from 
a standard argument in homological algebra. Here is a general version whose strength we will use later:

\begin{lemma} \label{cohomology chase} Let 
$$
\sM:\quad \cdots \to \sM^{-1} \to \sM^{0} \to \sM^{1} \to \cdots
 $$
be a complex of sheaves, and let $\sF^{i} := \rH^{i}\sM$ be the homology of $\sM$ at the $i$-th term. If 
\begin{align}
&H^{j-t} (\sM^{t})=0 \quad\text{ for all $t$,} \\
&H^{j- t-1}(\sF^{ t}) = 0 \;\text{ for all $t >0,$}\\
&H^{j-t+1}(\sF^t) = 0 \;\text{ for all $t <0$}
\end{align} then
 $H^j (\sF^{0})=0$.
\end{lemma}

%\frank{Should we mention $\dim \PP^n < \infty$ ?} \david{probably not!}

\begin{proof} Break $\sM$ into the short exact sequences
\begin{align*}
0\to &\sZ^{i} \to \sM^{i}\to\sB^{i+1}\to 0\\
0\to &\sB^{i}\to \sZ^{i}\to \sF^{i}\to 0
\end{align*}
and chase the corresponding long exact sequences in cohomology. 
\end{proof}

To complete the proof of Proposition~\ref{reg from monad}
we must show that if $\reg^{k}\sF\leq 0$ then
$\sF$ admits a linear monad with $\sM^{\ell} = 0$ for all $\ell>k$. The object we need is the one mentioned in \cite{EFS}  Example 8.5, and is constructed using the Beilinson-Gel'fand-Gel'fand correspondence (BGG). Since the property we need was not spelled out there, we review the construction and add some details.

Set $W=H^0 \sO_{\PP^{n}}(1)$, and let $E$ be the exterior algebra 
$E= \Lambda W^{*}$.
The cohomology table of a coherent sheaf $\sF$, as we have presented it, is also the
Betti table of the \emph{Tate resolution} $\TT(\sF) $, which is  a minimal graded free exact complex over $E$. The terms of $\TT(\sF)$ are
$$
\TT^e(\sF) = \oplus_{i=0}^n\Hom_{\mathbb K}(E,  H^i \sF(e-i))
$$
where $H^i \sF(e-i)$ is considered as a vector space concentrated in degree
$e-i$.
We consider the elements of $W^{*}$ as having degree $-1$, so the $E$-module
$\omega_{E} = \Hom_{\mathbb K}(E,K)$ is nonzero in degrees
$n+1, n\dots, 0$, and
$\Hom_{\mathbb K}(E,  H^i \sF(e-i))$
can be nonzero only in degrees $e-i+n+1,\dots, e-i$. 

Now suppose that $\reg^{k}\sF\leq 0$; this means that $H^{j}\sF(-j) =0$ for
$j> k$. Thus $\TT^0 \sF$ is generated in degrees $\geq -k+n+1$, and it follows
that the graded components of $\TT^0 \sF$ are all zero below degree $-k$. This
implies the same vanishing for the $E$-submodule
$P= \ker (\TT^0 \sF \to \TT^1 \sF)$. 

To the $E$-module $P$ the BGG correspondence associates a linear free complex 
$L(P)$ over $S$:
$$
L(P):  \cdots \to S\otimes P_{1}\rTo^{\partial} S\otimes P_{0}\rTo^{\partial} S\otimes P_{-1 }\rTo \cdots.
$$
The differential $\partial$ is defined to be multiplication by the element
$$
\sum_{i=0}^{n}x_{i}\otimes e_{i} \in S\otimes E,
$$
where $\{x_{i}\}$ and $\{e_{i}\}$ are dual bases of $W$ and $W^{*}$.
Since $P_{j}$ is concentrated in degree $j$, the module $S\otimes P_{j}$
is a direct sum of copies of $S(-j)$.

It follows from BGG (\cite{EFS} Theorem 8.1, with $\TT' = \TT(\sF)^{\geq 0}$ that 
the sheafification $\sM := \widetilde{L(P)}$ of $L(P)$ is a monad for $\sF$.
The term $\sM^{\ell}$ is equal to $\sO_{\PP^n}(\ell)^{\dim P_{-\ell}}$. 
The observation above that $P_{j} = 0$ for $j<-k$ implies that $\sM^{\ell} = 0$
for $\ell>k$, as required.
\end{proof}

The correspondence between the Tate resolution and the cohomology table allows us to generalize an important fact about Castelnuovo-Mumford regularity:

\begin{proposition}\label{inductive regularity}
If $\reg^{k}\sF = m$ then 
$H^{j} \sF(m'-j) = 0$  for all $j >k$ \emph{and} $m'\geq m$. Similarly, if
$\coreg^{\ell}\sF = m$ then
$H^{j} \sF(m-j) = 0$ for all $j <n-k$ \emph{and} $m'\leq m$.
\end{proposition}
\begin{proof} The given conditions with $m'=m$ are simply
the definitions of $\reg^{k}$ and $\coreg^{k}$.  If $H^{j} \sF(m'-j) \neq 0$ for some $m'>m$, then, because the Tate resolution is a minimal complex, no term
of the resolution could map into the summand  
$H:= \Hom_{\mathbb K}(E,  H^j \sF(m'-j))$, and it follows that this module would
be a submodule of one of the syzygies in the resolution. Since $H$
is an injective module over the exterior algebra, it would actually be a summand.
However, $H$ is also a free module over the exterior algebra, so this would contradict
the minimality of the Tate resolution. 

Since the dual of the Tate resolution is again exact and minimal, we can apply the
same argument to the dual to get the corresponding statement about coregularity.
\end{proof}

\section{Proof of Theorem~\ref{main2-intro} }\label{sec:main proofs}

The coregularity statement of Theorem~\ref{main2-intro} follows from the regularity
statement by dualising, so it is enough to prove the latter, in its strong
form. Thus we suppose that  $\sF$  and  $\sG$ are coherent sheaves on $\PP^n$
with $\dim Tor_{1}(\sF,\sG) \leq p+2$. It suffices to  show that
if $p = k+\ell$ then
$\reg^{p}(\sF\otimes \sG) \leq \reg^{k}\sF  + \reg^{\ell}\sG$.
Replacing $\sF$ and $\sG$ by $\sF(-k)$ and $\sG(-\ell)$ respectively, 
we may assume $\reg^{k}\sF = \reg^{\ell}\sG = 0$, and we must show that
for each $j>p$ we have
$H^{j}(\sF\otimes \sG(-j)) = 0$.

By Proposition \ref{reg from monad}, the sheaf $\sF$ has a linear monad 
of the form
$$
\sM: \cdots \to\sM^{-1} \to \sM^{0} \to \sM^{1} \to \cdots \to \sM^{k}\to 0
$$ 
where $\sM^t$ is a direct sum of copies of $\sO_{\PP^{n}}(t)$. 
Since the truncated complex 
$$
\sM^{+}:\quad \sM^{0} \to \cdots \to \sM^{k}\to 0
$$
is locally split, we have 
$$
\ker(\sM^{0}\to \sM^{1}) \otimes \sG = \ker(\sM^{0}\otimes \sG \to \sM^{1}\otimes \sG)
$$
and it follows that $\rH^{0}(\sM\otimes \sG(-j)) = \sF\otimes \sG(-j)$.

We now apply Lemma~\ref{cohomology chase} to 
the complex $\sM\otimes \sG(-j)$. 
Since the term $\sM^{t}\otimes \sG(-j)$ is a direct sum of copies of
$\sG(t-j)$,
 it suffices to show 
\begin{align}
&H^{j-t} \sG(t-j) = 0 \text{ for all $t\leq k$} \\
&H^{j-t-1}(\rH^{t} (\sM\otimes \sG(-j))) = 0 \text{ for all $t >0.$}\\
&H^{j-t+1}(\rH^{t} (\sM\otimes \sG(-j))) = 0 \text{ for all $t <0.$}
\end{align}
Since $j>p=k+\ell$ and $-t \ge -k$ we have $j-t > \ell$, and (4) holds because $\reg^\ell \sG=0$.

%To prove (3), note that $\reg^{\ell}\sG(t) = -t$, and it follows  that
%$\reg^{\ell'}\sG(t) \leq -t$ for all $\ell' \geq \ell.$ This amounts to saying that
%$$
%H^{\ell'-j}(\sG(t-j) = 0 \text{ for } j\geq 1 \text{ and } \ell'\geq \ell.
%$$
%Since $t\leq k$ we have $p-t \geq \ell$, and taking $\ell'=p-t$ completes the 
%proof of (3).

To prove (5) we observe that $\rH^{t}\sM\otimes \sG(-j) = 0$ for
all $t>0$ simply because $\sM^{+}$ is locally split. It remains to 
prove (6). But for $t<0$ and $j > p$ the number
$j-t+1\ge p+3$, so it is enough to show that
$\dim \rH^{t}(\sM\otimes \sG(-j)) = \dim \rH^{t}(\sM\otimes \sG) \leq p+2$.

The local splitting of $\sM^{+}$
further implies that $\sZ^{0} := \ker(\sM^{0}\to \sM^{1})$ is a vector
bundle, so 
$$
\sM^{-}: \quad \cdots \to \sM^{-2}\to\sM^{-1}\to \sZ^{0}
$$
is a locally free resolution of $\sF$. Thus for $t<0$
$$
\rH^{t}(\sM \otimes \sG) = Tor_{-t}(\sF, \sG)
$$
By the rigidity of Tor (Auslander, \cite{Auslander}), our hypothesis that
$\dim Tor_{1}(\sF, \sG)\leq p+2$ implies
$\dim Tor_{-t}(\sF, \sG)\leq p+2$
for all $t<0$, completing the proof.

\section{Proof of Theorem~\ref{sharpness}}

The statement for the coregularity follows from that for the regularity by duality, so we restrict ourselves to the regularity formulas.

We may shift $\sF$ and $\sG$ and assume without loss of generality that 
$\reg^{0}\sF = \reg^{0}\sG = 0$. By Theorem~\ref{BS for VB} we can write the
cohomology table $\Phi$ of $\sF$ as a sum of cohomology tables of 
homogeneous bundles in characteristic 0. Since
$H^{i}\sF(k) = 0$
for $i>0$ and $k\geq 0$, these bundles must in fact have $\reg^{0}\le 0$; that
is, they all have the form $S_{\lambda}Q$ for some partitions $\lambda$. Since $\reg^0 \sF =0$ we have at least one partition $\lambda$ with 
with $\lambda_{0} = 0$ occuring. Of course similar statements hold for the cohomology
table $\Gamma$ of $\sG$. 

Multiplying $\Phi$ and $\Gamma$ by sufficiently divisible integers, we may assume that the \BS decompositions have positive integral---not just rational---coefficients,
so that they correspond to actual homogeneous bundles. 

We will complete the proof of
Theorem~\ref{sharpness} 
by showing that if  $\sF$  is a direct sum of homogeneous bundles
\[
\sF  = \bigoplus_{u=0}^{v}S_{\lambda^{u}}, \text{ with } \lambda^{0}\leq \cdots\leq \lambda^{v}
\]
and similarly for $\sG$, then
$$
(*)\quad \reg^p(\sF\tensor \sG) = \min_{ k+l=p} \bigl(\reg^k\sF +\reg^\ell \sG \bigr)
$$
for every $0\leq p\leq n-1$. Since the inequality $\leq$ is part of
Theorem~\ref{main2-intro},
it suffices to show that the left hand side of $(*)$ is at least as large as the right hand side.

Since the $k$-regularity index of $S_{\lambda}Q$ is $-\lambda_{k}$, the minimum
on the right hand side of $(*)$ is achieved by the minimal partition 
involved in the decomposition. On the other hand, 
the $p$-th regularity index of a direct sum is the maximum of the $p$-th regularity
indices of the components, so it suffices to prove the inequality after replacing
each of $\sF$ and $\sG$ by a single summand, corresponding to the minimal
partitions in the two decompositions; that is, we may take
$\sF = S_{\lambda}Q$ and $\sG = S_{\mu}Q$ for some partitions $\lambda$ and
$\mu$.

We now have
$$
\sF\otimes \sG = S_{\lambda}Q\otimes S_{\mu}Q 
= \bigoplus S_{\nu^{u}}Q
$$
where the set of partitions $\nu^{u}$ (which may occur with multiplicity)
is determined by the Littlewood-Richardson formula. Since the 
regularity indices are the negatives of the parts of the partition, 
this translates into the following result in representation theory:

\begin{proposition}\label{rep theory}
Let $V$ be an $n$-dimensional vector space over a field of characteristic zero,
and let $0\leq p\leq n-1$. 
There is a representation $S_{\nu}V$ appearing in 
$S_{\lambda}V\otimes S_{\mu}V$, such that 
$\nu_{p}\leq \max_{ k+l=p} \lambda_{k}+\mu_{l}$.
\end{proposition}

\begin{proof}
Let
\begin{align*}
\lambda' &= \lambda_{p},\lambda_{p-1},\dots,\lambda_{0}\\
\mu' &= \mu_{p},\mu_{p-1},\dots,\mu_{0}.
\end{align*}
be the partitions obtained by truncating $\lambda$ and $\mu$.
One sees from the Littlewood-Richardson formula as described, for example,
in Fulton~\cite{Fulton}, that
if a representation corresponding to the partition $\nu'$ 
occurs in $S_{\lambda'}V\otimes S_{\mu'}V$
then the representation corresponding to the partition
$$
\nu = (\lambda_{n-1}+\mu_{n-1}, \dots, \lambda_{p+1}+\mu_{p+1}, \ 
\nu'_{p},\dots, \nu'_{0})
$$
occurs in $S_{\lambda}V\otimes S_{\mu}V$. Thus we may assume from
the outset that $p=n-1$.
If we set $g:= \max_{ k+l=n-1} \lambda_{k}+\mu_{l}$ then the termwise sum of
$\lambda$ with the sequence $(\mu_{0}, \dots, \mu_{n-1})$, which is
the reverse of $\mu$, is a sequence
of numbers $\leq g$. We want to show that  in the product
$S_{\lambda}V\otimes S_{\mu}V$ there occurs a representation
$S_{\nu}V$ such that $\nu_{n-1}\leq g$.

What we wish to prove can now be expressed as a statement about the intersection ring of
the Grassmannian $\Gr(n,n+g)$ of $n$-planes in $\bC^{n+g}$ as follows:
Let 
$$
\sV = (0\subsetneq V_{1}\subsetneq\cdots\subsetneq V_{n+g} = \bC^{n+g})
$$
be a complete flag in $\bC^{n+g}$.
We write $\Sigma_{\lambda}(\sV)$ for the Schubert cycle
in  $\Gr(n,n+g)$ defined by
$$
\Sigma_{\lambda}(\sV) = \{W\in \Gr(n, n+g)\mid
\dim W\cap V_{g+i-\lambda_{n-i}} \geq i \text{ for } 1\leq i\leq n-1\},
$$
and similarly for $\Sigma_{\mu}(\sV)$. The product of
the classes $[\Sigma_{\lambda}(\sV)]$ and 
$[\Sigma_{\mu}(\sV)]$ in the 
intersection ring of $\Gr(n,n+g)$ is the sum (with multiplicity) of 
the classes of those
$\Sigma_{\nu}$ such that $S_{\nu}V$ occurs in 
$S_{\lambda}V\otimes S_{\mu}V$ \emph{and} 
 $\nu_{n-1}\leq g$. (This is explained, and the proof sketched,
in \S~9 of \cite{Fulton}. Thus our problem is to show
that the
intersection product $[\Sigma_{\lambda}(\sV)][\Sigma_{\mu}(\sV)]$
is nonzero. This well-known fact can be proved as follows.

Choose another flag 
$$
\sV' = V_{0}'\subsetneq\cdots\subsetneq V_{n+g}'
$$ 
in general position with respect to $\sV$.
With the 
evident definition of $\Sigma_{\mu}(\sV')$, the product
above can
be computed as the class of the set-theoretic intersection
$$
\Sigma_{\lambda}(\sV)\cap \Sigma_{\mu}(\sV').
$$
(This follows, for example, from Kleiman's transversality theorem.)
Thus it suffices to show that this intersection is nonempty.

Since $\sV$ and $\sV'$ are in generic, the subspaces
$V_{i}\cap V_{n+g-i+1}'$ are all 1-dimensional. If
$e_{i}$ is a basis vector for this space,
then the conditions
$\lambda_{i}+\mu_{n-1-i}\leq g$ for $0\leq i\leq n-1$ imply
 that 
$$
W\in \Sigma_{\lambda}(\sV)\cap \Sigma_{\mu}(\sV'),
$$
where 
$$
W = \overline{e_{g+1-\lambda_{n-1}}, \dots, e_{g+n-\lambda_{0}}}
$$
so the product of the classes of these Schubert cycles
in nonzero as required.
\end{proof}

\section{Unobstructed Families of Vector Bundles}\label{sec:natural families}

Theorem \ref{main2-intro} gives a criterion for the vanishing of the obstruction space 
$\Ext^2(\sF,\sF)=0$.
In this section we describe the deformations of these unobstructed bundles. First the  criterion:

\begin{corollary}\label{Ext2} If $\sF$ be a vector bundle on $\PP^n$ with either
$\reg^0 \sF -\coreg^1 \sF \le 3$ \text{ or} $\reg^1 \sF -\coreg^0 \sF \le 3$, then
 the obstruction space $\rm Ext^2(\sF,\sF)$ vanishes.
\end{corollary}

\begin{proof} Since $\reg^k \sF^* = -\coreg^k \sF-1$  the assumption gives
$$\reg^1( \sF \tensor \sF^*) \le \min(\reg^0 \sF+\reg^1 \sF^*,\reg^1 \sF+\reg^0 \sF^*) \le 3-1=2$$
by Theorem \ref{main2-intro}.
Hence $\rm Ext^2(\sF,\sF)=H^2(\sF \tensor \sF^*) =0$. 
\end{proof}

Since replacing $\sF$ by $\sF^*$ interchanges the two assumption, we will focus on the case $\reg^1 \sF -\coreg^0 \sF \le 3$ in the following.
 To describe all bundles satisfying this assumption we will use Beilinson monad
\cite{EFS}, Theorem 6.1.
Given a sheaf $\sF$ on projective space the Beilinson monad
$$
\sB:  \quad \ldots \to \sB^{-1} \to \sB^0 \to \sB^1 \to \ldots
$$
 for $\sF$ has terms
$$
\sB^e = \oplus_j H^j(\sF(e-j)) \tensor \Omega^{j-e}(j-e).
$$
$\sB$ is obtained by applying the functor $\Omega$ to the Tate resolution $\TT(\sF)$,
where $\Omega$ is the additive functor that
sends the $E$-module $\omega_E(i)=\rm Hom_K(E,K(i))$ to the sheaf of twisted $i$-forms
$\Omega^i(i)$.
The identification 
$$
\rm Hom(\Omega^i(i),\Omega^j(j)) = \Lambda^{i-j} W^* = \rm Hom_E(\omega_E(i),\omega_E(j))
$$
provides the maps.

\begin{theorem}\label{monad of unobstructed}
Let $\sF$ be a vector bundle with $\reg^1 \sF -\coreg^0 \sF \le 3$ twisted such that $\reg^1 \sF=2$. Consider $A=\TT^{0} \sF$ and $B=\TT^{1} \sF$. There is non-empty Zariski open subset $U \subset \Hom_E(A,B)$ such that the kernel
$$  \sF_\varphi =\ker (\Omega(\varphi):\Omega A \to \Omega B) $$
is a vector bundle with the same Chern classes and rank as $\sF$.
\end{theorem}

\begin{proof} By assumption the cohomology of $\sF$ is non-zero in a range like the indicated one.
\medskip

\centerline{
\vbox{
\def\star{\rlap{\hbox to 13pt{\hfil\raise3.5pt\hbox{$*$}\hfil}}}
\def\ {\hbox to 13pt{\vbox to 13pt{}\hfil}}
\def\*{\star\ }
\def\a{{\rlap{\hbox to 13pt{\hfil\raise3.5pt\hbox{$a$}\hfil}}}\ }
\def\_{\hbox to 13pt{\hskip-.2pt\vrule\hss\vbox to 13pt{\vskip-.2pt
            \hrule width 13.4pt\vfill\hrule\vskip-.2pt}\hss\vrule\hskip-.2pt}}
\def\x{\star\_}
\offinterlineskip
\hbox{\_\_\_\_\x\x\ \ \ \  \  \ \ \ \ \ \ \ \ \ \ \ }
\hbox{\ \ \ \ \x\x  \\ \ \   \\ \ \   \ \ \ \ \ \ \ \ \ \ \ \  }
\hbox{\ \ \ \ \x \x  \\ \ \   \\ \ \   \ \ \ \ \ \ \ \ \ \ \ \ }
\hbox{\ \ \ \ \x \x  \\ \ \   \\ \ \   \ \ \ \ \ \ \ \ \ \ \ }
\hbox{\ \ \ \ \x \x  \_\_\_ \   \\ \ \   \ \ \ \ \ \ \ \ \ \ }
\hbox{\ \ \ \  \x\_\_\_\_\_\_\_\_\_   \   \  \  \  \ }
}
}  
\noindent 
The Beilinson monad of $\sF$ depends only the terms in the range indicated by a $*$,  since all other terms are zero or are sent to zero. In particular $\Omega(\TT(\sF))$ is a two term complex
$$ 
0 \to \Omega A \to \Omega B \to 0
$$
Since it is an open condition for $\varphi$ to define a monad
$ 0 \to \Omega A \to \Omega B \to 0$
and the set of $\varphi$ is non empty by the existence of $\sF$, the result follows.
\end{proof}

Recall from \cite{ES1} that a vector bundle $\sF$ on $\PP^{n}$ has \emph{natural cohomology} in the sense of Hartshorne and Hirschowitz if, 
for each $d\in \ZZ$, at most one of the groups $H^i\sF(d)$ is non-zero. The bundles $\sF$ is called  \emph{supernatural} if, in addition, the polynomial function $\chi(\sF(d))$ has $n$ distinct integral roots.

\begin{corollary}\label{natural unobstructed} Let $\sF$ be a vector bundle with $\reg^1 \sF -\coreg^0 \sF \le 3$ normalized (by tensoring with a line bundle)
so that $\reg^1 \sF=2$, and assume that $\sF$ has natural cohomology. 
Every vector bundle with natural cohomology with the same rank,  Chern classes, regularity and coregularity indices arises as 
$$  \sF_\varphi =\ker (\Omega(\varphi):\Omega A \to \Omega B) $$
for some $\varphi \in U$. In particular these vector bundles form an irreducible unirational family.
\end{corollary}

\begin{proof} The cohomology table of any of these bundles is determined by the Hilbert polynomial $d \mapsto \chi(\sF(d))$, since for each twist at most two terms could be nonzero due to the narrow cohomology river, and because we have natural cohomology. So they all have the same cohomology table and all arise from some $\varphi \in U$.
\end{proof}

Note that these bundles are not necessarily stable. For example we could have a direct sums of bundles with different slopes in these family. Thus we do not speak of a moduli space.

Corollary \ref{natural unobstructed}  does not settle the existence problem for such bundles. However, it provides a computational criterion:
A bundle with desired unobstructed natural cohomology table exists if and only if a general map
$\varphi \in U$ yields such a bundle. \BS theory characterizes the cohomology tables that can occur, up to a rational multiple. Given an integral cohomology table $\gamma$
satisfying the numerical condition of Corollary~\ref{natural unobstructed} such that some
multiple of $\gamma$ is the cohomology table of a bundle, we conjecture that there is a number $c_{0}(\gamma)$ such that  $c \gamma$ is the cohomology table of a bundle if and only if $c \ge c_0(\gamma)$, as in the following example.

\begin{example} The table $\gamma$
\begin{verbatim}
       4: 56 15 . . .  .
       3:  .  . 2 . .  .
       2:  .  . . 1 .  .
       1:  .  . . . .  .
       0:  .  . . . 8 35
          -2 -1 0 1 2  3
\end{verbatim}
``looks like'' the cohomology table of a rank 4 vector bundle on $\PP^4$, but
it is not! This is because
for any two  2-forms $(\eta_1,\eta_2) \in \oplus_1^2 \Lambda^2 W^* \subset E^2$
the kernel of the wedge product
$\ker( \Lambda^2  W^* \to  \oplus_1^2 \Lambda^4 W^*)$ is always nonzero.

Indeed, $\eta_i \in \Lambda^2 V_i$ for some 4 dimensional subspace $V_i \subset W^*$
and the annihilator of $\eta_i$ has codimension one in $\Lambda^2 V_i$.
Thus the intersection $\ann(\eta_i) \cap \Lambda^2 (V_1 \cap V_2)$ has codimension at most 1
in the 3-dimensional space $\Lambda^2 (V_1 \cap V_2)$ and the intersection $
\ann(\eta_1) \cap \ann(\eta_2)$ is at least 1-dimensional.

However, experiments with Macaulay2~\cite{M2} convince us that every multiple $c\gamma$ with $c \ge 2$ does occur as the cohomology table of a bundle.\end{example}

\bigskip

\vbox{\noindent Author Addresses:\par
\smallskip
\noindent{David Eisenbud}\par
\noindent{Department of Mathematics, University of California, Berkeley,
Berkeley CA 94720}\par
\noindent{eisenbud@math.berkeley.edu}\par
\smallskip
\noindent{Frank-Olaf Schreyer}\par
\noindent{Mathematik und Informatik, Universit\"at des Saarlandes, Campus E2 4, 
D-66123 Saarbr\"ucken, Germany}\par
\noindent{schreyer@math.uni-sb.de}\par
}

\end{document}